\documentclass[12pt,leqno]{amsart}
\usepackage{oldgerm}
\usepackage{latexsym}
\usepackage{amsmath}
\usepackage{amssymb}
\usepackage{delarray}
\usepackage{amscd}
\usepackage{amsthm}

\begin{document} 
 \title[Congruence]{Congruence between Duke-Imamoglu-Ikeda  lifts and non-Duke-Imamoglu-Ikeda  lifts} 
\author{Hidenori Katsurada} 

\date{January 10, 2011}

\maketitle

\def\vecttwo(#1;#2){\left(\begin{array} {c} #1 \\ #2
\end{array}\right)}

\def\vecttwobrack(#1;#2){\left[\left(\begin{array} {c} #1 \\ #2
\end{array}\right)\right]}

\def\vectthree(#1;#2;#3){\left(\begin{array} {c} #1 \\ #2 \\ #3
\end{array}\right)}

\def\vectfour(#1;#2;#3;#4){\left(\begin{array} {c} #1 \\ #2 \\ #3 \\
#4 \end{array}\right)}

\def\mattwo(#1;#2;#3;#4){\left(\begin{array}{cc} #1 & #2 \\ #3 &#4
\end{array}\right)}

\def\matthree(#1;#2;#3;#4;#5;#6;#7;#8;#9){\left(\begin{array}{ccc} 
#1 & #2 & #3 \\ 
#4 &#5  & #6 \\
#7 &#8 & #9
\end{array}\right)}

\def\diag{{\rm diag}}

\bigskip
   
\begin{abstract} Let $k$ and $n$ be positive even integers. For a Hecke eigenform $g$ in the Kohnen plus sbspace of weight $k-n/2+1/2$ for $\varGamma_0(4),$  let $f$ be the primitive form of weight $2k-n$ for $SL_2({\bf Z})$ corresponding to $g$ under the Shimura correspondence, and $I_n(g)$ the Duke-Imamoglu-Ikeda lift of $g$ to the space of cusp forms of weight $k$ for $Sp_n({\bf Z}).$  Then we characterize prime ideals giving congruence between $I_n(g)$ and another cuspidal Hecke eigenform not coming from the Duke-Imamoglu-Ikeda lift in terms of the  specilal values of the  Hecke $L$-function and the adjoint $L$-function of $f.$
\end{abstract}

\bigskip

\section{Introduction}

As is well known, there is a congruence between the Fourier coefficients of Eisenstein series and those of cuspidal Hecke eigenforms for $SL_2({\bf Z}).$ This type of congruence is not only interesting in its own right but also has an importan application to number theory as shown by Ribet [Ri1]. Since then, there have been so many important results about the congruence of elliptic modular forms (cf. [Hi1],[Hi2],[Hi3].) In the case of Hilbert modular forms or Siegel modular forms, it is sometimes more natual and important to consider the congruence between the Hecke eigenvalues of Hecke eigenforms modulo a prime ideal $\textfrak P$. We call such a $\textfrak P$ a prime ideal giving the congruence or a congruence prime. 
For a cuspidal Hecke eigenform $f$ for $SL_2({\bf Z}),$ let $\hat f$ be a lift of $f$ to the space ${\textfrak M}_l(\varGamma')$ 
of modular forms of weight $l$ for a modular group $\varGamma'.$ Here we mean by the lift of $f$ a cuspidal Hecke eigenform whose ceratin L-function can be expressed in terms of certain L-functions of $f.$ As typical examples of the lift we can consider the Doi-Naganuma lift,  the Saito-Kurokawa lift, and the Duke-Imamoglu-Ikeda  lift. We then consider the following problem:

\bigskip

{\bf Problem.} {\it Characterize the prime ideals giving the congruence between $\hat f$ and a cuspidal Hecke eigenform in  ${\textfrak M}_l(\varGamma')$ not coming from the lift. In particular characterize them in terms of special values of certain L-functions of $f.$}

\bigskip

This type of problem was first invetigated in the Doi-Naganuma lift case by Doi, Hida, and Ishii [D-H-I]. 
In our previous paper [Ka2], we considered the relationship between the congruence of cuspidal Hecke eigenforms with respect to $Sp_n({\bf Z})$ and the special values of their standard zeta functions.  In particular, we proved a conjecture  proposed by Harder [Ha] concerning the congruence between Saito-Kurokawa lifts and non-Saito-Kurokawa lifts under certain conditions. (See also [Br].) In this paper,  we consider  a congruence between Duke-Imamoglu-Ikeda  lifts and non-Duke-Imamoglu-Ikeda  lifts, which is a generalization of our previous conjecture. We explain our main result more precisely. Let $k$ and $n$ be positive even integers. For a Hecke eigenform $g$ in the Kohnen plus sbspace of weight $k-n/2+1/2$ for $\varGamma_0(4),$  let $f$ be the primitive form of weight $2k-n$ for $SL_2({\bf Z})$ corresponding to $g$ under the Shimura correspondence, and $I_n(g)$ the Duke-Imamoglu-Ikeda lift of $g$ to the space of cusp forms of weight $k$ for $Sp_n({\bf Z}).$  Moreover let $L(s,f)$ and $L(s,f,{\rm Ad})$ be the Hecke's $L$-function and the adjoint $L$-function of $f,$ respectively. Then our main result can be roughly stated as follows:

\bigskip

{\it A prime ideal in the Hecke field of $f$ dividing the numerator of the algebraic part of $L(k,f)\prod_{i=1}^{n/2-1}L(2i+1,f,{\rm Ad})$ gives a congruence between $I_n(g)$ and a Hecke eigenform not coming from the Duke-Imamoglu-Ikeda lift }(cf. Theorem 4.7).

\bigskip
 The paper is organized as follows. In Section 2, we review the Hecke theory of Siegel modular forms and the standand $L$ function.  In Section 3, we review a result concerning the relationship between the congruence of cuspidal Hecke eigenforms with respect to $Sp_n({\bf Z})$ and the special values of their standard zeta functions. In Section 4, we propose a conjecture concerning the congruence between  Duke-Imamoglu-Ikeda  lifts and non-Duke-Imamoglu-Ikeda  lifts, and prove it under a certain condition. By using this theorem, we can give a non-trivial example of congruence between the Duke-Imamoglu-Ikeda lift and non-Duke-Imamoglu-Ikeda lift (cf. Example in Section 4). 
 
The author thanks Professor C. Poor and Professor D. Yuen for informing him the estimate of the dimension of  the space of  cusp forms of weight $18$ for $Sp_4({\bf Z}).$   The author also thanks Professor B. Heim, Professor H. Hida, Professor S. Yasuda, and Professor T. Yamauchi for their valuable comments.

{\bf Notation.}   For a
commutative ring $R$, we denote by $M_{mn}(R)$ the set of
$(m,n)$-matrices with entries in $R.$  In particular put $M_n(R)=M_{nn}(R).$ 
Here we understand
$M_{mn}(R)$ the set of the {\it empty matrix} if $m=0$ or
$n=0.$  For an $(m,n)$-matrix $X$ and an $(m,m)$-matrix
$A$, we write $A[X] = {}^t X A X,$ where $^t X$ denotes the
transpose of $X$.  Let $a$ be an element of $R.$ Then
for an element $X$ of $M_{mn}(R)$ we often use the same
symbol $X$ to denote the coset $X \ {\rm mod} \ aM_{mn}(R).$ Put $GL_m(R) = \{A \in M_m(R)  \ | \  \det A \in R^* \}$, where $\det
A$ denotes the determinant of a square matrix $A$, and $R^*$
denotes the unit group of $R.$  
 Let $S_n(R)$ denote
the set of symmetric matrices of degree $n$ with entries in
$R.$ Furthermore, for an integral domain $R$ of characteristic different 
from $2,$  let ${\mathcal H}_n(R)$
denote the set of half-integral matrices of degree $n$ over
$R,$ that is, ${\mathcal H}_n(R)$ is the set of symmetric
matrices of degree $n$ whose $(i,j)$-component belongs to
$R$ or ${1 \over 2}R$ according as $i=j$ or not.  For a subset $S$ of
$M_n(R)$ we denote by $S^{\times}$ the subset of $S$
consisting of non-degenerate matrices.  In particular, if
$S$ is a subset of $S_n({\bf R})$ with ${\bf R}$ the field
of real numbers, we denote by $S_{>0}$ (resp. $S_{\ge 0}$) the subset of $S$
consisting of positive definite (resp. semi-positive definite) matrices.  Let $R'$ be a
subring of $R.$ Two symmetric matrices $A$ and $A'$ with
entries in $R$ are called equivalent over $R'$ with each
other and write $A {\sim}_{R'} A'$ if there is
an element $X$ of $GL_n(R')$ such that $A'=A[X].$ We also write $A \sim A'$ if there is no fear of confusion. 
For square matrices $X$ and $Y$ we write $X \bot Y = \mattwo(X;O;O;Y).$

\section{Standard zeta functions of Siegel modular forms} 
For a complex number $x$ put ${\bf e}(x)=\exp(2\pi \sqrt{-1} x).$ Furthermore put $J_n=\mattwo(O_n;-1_n;1_n;O_n),$ where $1_n$ denotes the unit matrix of degree $n.$ For a subring $K$ of ${\bf R}$ put  
$$GSp_n(K)^+=\{M \in GL_{2n}(K)  \ | \  J_n[M]= \kappa(M)J_n \
{\rm with \ some} \ \kappa(M) > 0 \},$$
 and
$$Sp_n(K)=\{M \in GSp_n(K)^+  \ | \  J_n[M]=J_n \}. $$
 Furthermore, put 
$$\varGamma^{(n)}=Sp_n({\bf Z})=\{M \in GL_{2n}({\bf Z})   \ | \  J_n[M]=J_n \}. $$
We sometimes wrtite an element $M$ of $GSp_n(K)$ as 
$M=\mattwo(A;B;C;D)$ with $A,B,C,D \in M_2(K).$ We define a subgroup $\varGamma^{(n)}_0(N)$ of $\varGamma^{(n)}$ as
$$\varGamma^{(n)}_0(N)=\{\mattwo(A;B;C;D) \in \varGamma^{(n)} \ | \ C \equiv O_n \ {\rm mod} \ N \}.$$
 Let ${\bf H}_n$ be Siegel's
upper half-space. For each element $M =\mattwo(A;B;C;D) \in GSp_n({\bf R})^+$
and $Z \in {\bf H}_n$ put 
$$M(Z)=(AZ+B)(CZ+D)^{-1}$$ and 
$$j(M,Z)=\det (CZ+D).$$
Furthermore, for a function $f$ on ${\bf H}_n$ and an integer $k$ we define $f|_k M$ as
$$(f|_k M)(Z)= \det (M)^{k/2} j(M,Z)^{-k}f(M(Z)).$$
For an integer or half integr $l$ and the subgroup $\varGamma^{(n)}_0(N)$ of $\varGamma^{(n)},$ we denote by ${\textfrak M}_{k}(\varGamma^{(n)}_0(N))$ the space of (holomorphic)  modular forms of weight $k$
 with respect to $\varGamma^{(n)}_0(N).$   We denote by ${\textfrak S}_{k}(\varGamma^{(n)}_0(N))$ the sub-space of ${\textfrak M}_{k}(\varGamma^{(n)}_0(N))$ consisting of cusp forms. 
Let $f$ be a modular form  of weight $k$ with respect to  $\varGamma^{(n)}_0(N).$ Then $f$ has the following Fourier expansion: 
$$f(Z)=\sum_{A \in {\mathcal H}_n({\bf Z})_{\ge 0}} c_f(A) {\bf e}({\rm tr}(AZ)),$$ 
and in particular if $f$ is a cusp form, $f$ has the following Fourier expansion:
$$f(Z)=\sum_{A \in {\mathcal H}_n({\bf Z})_{> 0}} c_f(A) {\bf e}({\rm tr}(AZ)),$$ where ${\rm tr}$ denotes the trace of a matrix.
Let $dv$ denote the invariant volume element on ${\bf H}_n$ defined by 
$$dv=\det ({\rm Im}(Z))^{-n-1} \wedge_{1 \le j \le l \le n}(dx_{jl} \wedge dy_{jl}).$$
 Here for $Z \in {\bf H}_n$ we write $Z= (x_{jl}) + \sqrt{-1} (y_{jl})$ with real matrices $(x_{jl})$ and $(y_{jl}).$ For two modular forms $f$ and $g$ of weight $l$ with respect to $\varGamma^{(n)}_0(N)$ we define the Petersson scalar product $\langle f,g \rangle $ by 
$$\langle f,g \rangle =[\varGamma^{(n)}:\varGamma^{(n)}_0(N)]^{-1}\int_{\varGamma^{(n)}_0(N) \backslash {\bf H}_n} f(Z)\overline {g(Z)} \det ({\rm Im}(Z))^l dv,$$
provided the integral converges.

Let ${\bf L}_n={\bf L}_{\bf Q}(GSp_n({\bf Q})^+,\varGamma^{(n)})$
denote the Hecke algebra over ${\bf Q}$ associated with the Hecke 
pair $(GSp_n({\bf Q})^+,\varGamma^{(n)}).$  Furthermore, let ${\bf L}_n^{\circ}={\bf L}_{\bf Q}(Sp_n({\bf Q}),\varGamma^{(n)})$
denote the Hecke algebra over ${\bf Q}$ associated with the Hecke 
pair $(Sp_n({\bf Q}),\varGamma^{(n)}).$ For
each integer $m$ define an element $T(m)$ of ${\bf L}_n$
by 
$$T(m)=\sum_{d_1,...,d_n,e_1,...,e_n}\varGamma^{(n)}(d_1 \bot ...\bot d_n \bot e_1 \bot ...\bot e_n)\varGamma^{(n)},$$ where
$d_1,...,d_n,e_1,...,e_n$ run over all positive integer satisfying
$$d_i|d_{i+1}, \ e_{i+1}|e_i \ (i=1,...,n-1), d_n|e_n,d_ie_i=m \
(i=1,....,n).$$ Furthermore, for $i=1,...,n$ and a prime number $p$
 put
$$T_i(p^2)=\varGamma^{(n)}(1_{n-i} \bot p1_i \bot p^21_{n-i} \bot p1_i)\varGamma^{(n)},$$
and $(p^{-1})=\varGamma^{(n)}(p^{-1} 1_{n})\varGamma^{(n)}.$
As is well known, ${\bf L}_n$ is generated over ${\bf Q}$ by all $T(p),
T_i(p^2) \ (i=1,...,n),$ and $(p^{-1}).$ We denote by ${\bf L}_n'$ the subalgebra of ${\bf L}_n$ generated by over ${\bf Z}$ by all $T(p)$ and
$T_i(p^2) \ (i=1,...,n).$ 
Let $T =\varGamma^{(n)} M \varGamma^{(n)}$ be an element of ${\bf L}_n \otimes {\bf C}.$ Write $T$ as $T=\cup_{\gamma} \varGamma^{(n)} \gamma$ and for $f \in {\textfrak M}_k(\varGamma^{(n)})$ define the Hecke operator $|_k T$ associated to $T$ as
$$f|_kT =\det (M)^{k/2-(n+1)/2}\sum_{\gamma} f|_k \gamma.$$ 
 We call this action  the Hecke operator as usual (cf. [A].) 
 If $f$ is an eigenfunction of a Hecke operator $T \in {\bf L}_n \otimes {\bf C}$,  we denote by $\lambda_f(T)$ its eigenvalue. Let ${\bf L}$ be a subalgebra of ${\bf L}_n.$ We call $f \in {\textfrak M}_k(\varGamma^{(n)})$ a Hecke eigenform for ${\bf L}$ if it is a common eigenfunction of all Hecke operators in ${\bf L}.$ In particular if ${\bf L}={\bf L}_n$ we simply call $f$ a Hecke eigenform.
  Furthermore, we denote by ${\bf Q}(f)$ the field generated over ${\bf Q}$ by eigenvalues of all $T \in {\bf L}_n$ as in Section 1.
 As is well known, ${\bf Q}(f)$ is a totally real algebraic number field of finite degree.  Now, first we consider the integrality of the eigenvalues of Hecke operators. For an algebraic number field $K,$ let ${\textfrak O}_K$ denote the ring of integers in $K.$ The following assertion has been proved in  [Mi2] (see also [Ka2].)  
 
 \bigskip
 
 {\bf Theorem 2.1} {\it Let $k \ge n+1.$ Let $f \in  \textfrak{S}_{k}(\varGamma^{(n)})$ be a common eigenform in ${\bf L}_n'.$ Then $\lambda_f(T)$ belongs to ${\textfrak O}_{{\bf Q}(f)}$ for any $T \in {\bf L}_n'.$ }

\bigskip

   Let ${\bf
L}_{np}={\bf L}(GSp_n({\bf Q})^+  \cap GL_{2n}({\bf Z}[p^{-1}]), \varGamma^{(n)})$ be the Hecke algebra associated with the pair $(GSp_n({\bf Q})^+  \cap GL_{2n}({\bf Z}[p^{-1}]), \varGamma^{(n)}).$  ${\bf
L}_{np}$ can be considered as a subalgebra of ${\bf L}_{n},$ and is generated over ${\bf Q}$ by $T(p)$ and $T_i(p^2) \ (i=1,2,.., n),$ and $(p^{-1}).$ We now review the Satake $p$-parameters of ${\bf L}_{np};$ let ${\bf P}_n={\bf Q}[X_0^{\pm},X_1^{\pm},...,X_n^{\pm}]$ be the ring of Laurent polynomials in $X_0,X_1,...,X_n$ over ${\bf Q}.$ Let ${\bf W}_n$ be the group of ${\bf Q}$-automorphisms of ${\bf P}_n$ generated by all permutations in variables $X_1,....,X_n$ and by the automorphisms $\tau_1,....,\tau_n$ defined by
$$\tau_i(X_0)=X_0X_i,\tau_i(X_i)=X_i^{-1},\tau_i(X_j)=X_j \ (j\not=i).$$
Furthermore, a group $\tilde {\bf W}_n$ isomorphic to ${\bf W}_n$ acts on the set $T_n=({\bf C}^{\times})^{n+1}$ in a way similarly to above. 
Then there exists a  ${\bf Q}$-algebra isomorphism $\Phi_{np}$, called the Satake isomorphism, from  ${\bf L}_{np}$ to the ${\bf W}_n$-invariant subring ${\bf P}_n^{{\bf W}_n}$ of ${\bf P}_n.$  Then for a ${\bf Q}$-algebra homomorphism $\lambda$ from ${\bf L}_{np}$ to ${\bf C},$   there exists an element $(\alpha_0(p,\lambda),\alpha_1(p,\lambda),...,\alpha_n(p,\lambda))$ of ${\bf T}_n$ satisfying
$$\lambda(\Phi_{np}^{-1}(F(X_0,X_1,...,X_n)))=F(\alpha_0(p,\lambda),\alpha_1(p,\lambda),...,\alpha_n(p,\lambda))$$
for $F \in  {\bf P}_n^{{\bf W}_n}.$ The equivalence class of 
$(\alpha_0(p,\lambda),\alpha_1(p,\lambda),...,\alpha_n(p,\lambda))$ under the action of $\tilde {\bf W}_n$ is uniquely determined  by $\lambda.$ We call this the Satake parameters of ${\bf L}_{np}$ determined by $\lambda.$

 Now
assume that an element $f$ of $M_k(Sp_n({\bf Z}))$ is a Hecke eigenform. Then for each prime number $p,$ $f$ defines a ${\bf Q}$-algebra homomorphism $\lambda_{f,p}$ from ${\bf L}_{np}$ to ${\bf C}$ in a usual way, and we denonte by $\alpha_0(p),\alpha_1(p),....,\alpha_n(p)$ the Satake parameters of ${\bf L}_{np}$ determined by
$f.$ We then define the standard zeta function $L(s,f,\underline {\rm St})$
by  $$L(s,f,\underline {\rm St})=\prod_p \prod_{i=1}^n
\{(1-p^{-s})(1-\alpha_i(p)p^{-s})(1-\alpha_i(p)^{-1}p^{-s})\}^{-1}.$$
 Let $f(z)=\sum\limits_{A \in {\mathcal H}_n({\bf Z})_{>0}} a(A){\bf e}({\rm tr}(Az))$ be a Hecke eigenform in $\textfrak{S}_{k}(\varGamma^{(n)}).$  
For a positive integer $m \le k-n$ such that $m \equiv n \ {\rm mod} \ 2$ put 
$$\Lambda(m,f,\underline {\rm St})= (-1)^{n(m+1)/2+1}2^{-4kn+3n^2+n+(n-1)m+2}$$
$$ \times \Gamma(m+1)\prod_{i=1}^n\Gamma(2k-n-i) {L(m,f,\underline {\rm St}) \over <f,f>\pi^{-n(n+1)/2+nk+(n+1)m}}.$$
Then the following theorem is due to B{\"o}cherer [B2] and Mizumoto [Mi].
  
\bigskip

{\bf Theorem  2.2.} {\it Let $l,k$ and $n$ be a positive integers  such that $\rho(n) \le l \le k-n,$ where $\rho(n)=3,$ or 1 according as $n \equiv 1 \ {\rm mod} \ 4$ and $n \ge 5,$ or not. 
Let $f \in \textfrak{S}_{k}(\varGamma^{(n)})$ be a Hecke eigenform. Assume that all the Fourier coefficients of  $f$ belong to ${\bf Q}(f).$  Then  $\Lambda(m,f,\underline {\rm St})$ belongs to ${\bf Q}(f).$
}

\bigskip

For later purpose, we consider a special element in ${\bf L}_{np};$ the polynomial $X_0^2X_1X_2 \cdots X_n \sum_{i=1}^n (X_i+X_i^{-1})$ is an element of ${\bf P}_n^{{\bf W}_n},$ and thus we can define an element $\Phi_{np}^{-1}(X_0^2X_1X_2 \cdots X_n \sum_{i=1}^n (X_i+X_i^{-1}))$ of ${\bf L}_{np},$ which is denoted by ${\bf r}_1.$

\bigskip

{\bf Proposition 2.3.} {\it Under the above notation the element ${\bf r}_1$ belongs to ${\bf L}_n',$ and we have
$$\lambda_f({\bf r}_1)=p^{n(k-(n+1)/2)}\sum_{i=1}^n (\alpha_i(p)+\alpha_i(p)^{-1}).$$}

\bigskip

{\it Proof.}  By a careful analysis of the computation in page 159-160 of [A], we see that ${\bf r}_1$ is a ${\bf Z}$-linear combination of $T_i(p^2) \ (i=1,...,n),$ and therefore we can prove the first assertion. Furthermore, by Lemma 3.3.34 of [A], we can prove the second assertion.

\section{Congruence of modular forms and special values of the standard zeta functions}

In this section we review a result concerning the congruence between the Hecke eigenvalues of modular forms of the same weight following [Ka2]. Let $K$ be an algebraic number field, and ${\textfrak O}={\textfrak O}_K$  the ring of integers in $K.$ For  a prime ideal ${\textfrak P}$ of ${\textfrak O},$ we denote by ${\textfrak O}_{({\textfrak P})}$ the localization of ${\textfrak O}$ at ${\textfrak P}$ in $K.$ Let ${\textfrak A}$ be a fractional ideal in $K.$ If ${\textfrak A}={\textfrak P}^e {\textfrak B}$ with ${\textfrak B}{\textfrak O}_{({\textfrak P})}={\textfrak O}_{({\textfrak P})}$  we write ${\rm ord}_{{\textfrak P}}=e.$ We simply write ${\rm ord}_{{\textfrak P}}(c)= {\rm ord}_{{\textfrak P}}((c))$ for $c \in K.$  
Now let $f$ be a Hecke eigenform in ${\textfrak S}_k(\varGamma^{(n)})$ and $M$ be a subspace of 
 ${\textfrak S}_k(\varGamma^{(n)})$ stable under Hecke operators $T \in {\bf L}_n.$ Assume that $M$ is contained in $({\bf C}f)^{\bot},$ where $({\bf C}f)^{\bot}$
  is the orthogonal complement of ${\bf C}f$ in ${\textfrak S}_k(\varGamma^{(n)})$ with respect to the Petersson product. Let $K$ be an algebraic number field containing ${\bf Q}(f).$ A prime ideal ${\textfrak P}$ of ${\textfrak O}_{K}$ is called a congruence prime of $f$ with respect to $M$ if there exists a Hecke eigenform $g \in M$   such that 
 $$\lambda_f(T) \equiv \lambda_g(T) \ {\rm mod} \ \tilde {\textfrak P}$$
 for any $T \in {\bf L}_n',$ where $\tilde {\textfrak P}$ is the prime ideal of
 ${\textfrak O}_{K {\bf Q}(g)}$ lying above ${\textfrak P}.$ 
 If $M= ({\bf C}f)^{\bot},$ we simply call ${\textfrak P}$ a congruence prime of $f.$ 
 
Now we consider the relation between the congruence primes and the standard zeta values. To consider this, we have to normalize the standard zeta value  $\Lambda(l,f,\underline {\rm St})$ for a Hecke eigenform $f$ because it is not uniquely determined by the system of Hecke eigenvalues of $f.$ We note that  there is no reasonable normalization of cuspidal Hecke eigenform in the higher degree case unlike the elliptic modular case.  Thus we define the following quantities: for a  Hecke eigenform $f(z)=\sum_{A} c_f(A){\bf e}({\rm tr}(Az))$ in ${\textfrak S}_{k}(\varGamma^{(n)}),$  let ${\textfrak I}_f$ be the ${\textfrak O}_{{\bf Q}(f)}$-module  generated by all $c_f(A)'s.$  Assume that the following:

(*) {\it There exists  a complex number $c$ such that all the Fourier coefficients $cf$ belongs to ${\bf Q}(f).$}

 Then ${\textfrak I}_f$ is a fractional ideal in ${\bf Q}(f),$ and therefore, so is $\Lambda(l,f,\underline {\rm St}){\textfrak I}_f^2$ if $l$ satisfies the condition in Theorem 2.2. We note that this fractional ideal does not depend on the choice of $c.$ We also note that the value
 $N_{{\bf Q}(f)}(\Lambda(l,f,\underline {\rm St}))N({\textfrak I}_f)^2$ does not depend on the  choice of $c,$ where $N({\textfrak I}_f)$ is the norm of the ideal ${\textfrak I}_f.$  Then, by [Kat2],  we have 
  
\bigskip
 
{\bf Theorem 3.1.} {\it  Let $f$ be a Hecke eigenform in ${\textfrak S}_k(\varGamma^{(n)}).$ Assume the condition (*).  Let $l$ be a positive integer satisfying the condition in Theorem 2.2. Let ${\textfrak P}$ be a prime ideal of ${\textfrak O}.$ Assume that ${\rm ord}_{{\textfrak P}}(\Lambda(l,f,\underline {\rm St}){\textfrak I}_f^2)<0$ and that it does not divide $(2l-1)!.$ Then ${\textfrak P}$ is a congruence prime of $f.$ In particular, if a rational prime number $p$ divides the denominator of $N_{{\bf Q}(f)}(\Lambda(l,f,\underline {\rm St}))N({\textfrak I}_f)^2,$ then $p$ is divisible by some congruence prime of $f.$ } 
 
\bigskip

{\bf Remark.} If the multiplicity one property holds for $S_k(\varGamma^{(n)}),$ the condition (*) holds for any cuspidal Hecke eigenform in $S_k(\varGamma^{(n)}).$ In [Kat3], we have treated the case where the condition (*) does not necessarily  hold.

\bigskip 

Now for a Hecke eigenform $f$ in ${\textfrak S}_k(\varGamma^{(n)}),$ let ${\textfrak T}_f$ denote the subspace of ${\textfrak S}_k(\varGamma^{(n)})$ spanned by all  Hecke eigenforms with the same system of the Hecke eigenvalues as $f.$

{\bf Corollary.} {\it In addition to the above notation and the assumption, assume that ${\textfrak S}_k(\varGamma^{(n)})$ has the multiplicity one property. Then ${\textfrak P}$ is a congruence prime of $f$ with respect to ${\textfrak T}_f^{\bot}.$ In particular, if a rational prime number $p$ divides the denominator of $N_{{\bf Q}(f)}(\Lambda(l,f,\underline {\rm St}))N({\textfrak I}_f)^2,$ then $p$ is 
divisible by some congruence prime of $f$ with respect to ${\textfrak T}_f^{\bot}$} 
 
\bigskip

\section{Congruence between Duke-Imamoglu-Ikeda  lifts and non-Duke-Imamoglu-Ikeda  lifts}
In this section, we consider the conguence between Duke-Imamoglu-Ikeda  lifts and non-Duke-Imamoglu-Ikeda  lifts. Throughout this section and the next, we assume that $n$ and $k$ are even positive integers.
Let $g$ be a Hecke eigenform belonging to the Kohnen plus space ${\textfrak S}_{k-n/2+1/2}^+(\varGamma_0(4)).$ Then $g$  has the following Fourier expansion:  
$$g(z)=\sum_{e}c_g(e){\bf e}(ez),$$ 
where $e$ runs over all positive integers such that $(-1)^{k-n/2}e \equiv 0, 1 \ {\rm mod} \ 4.$ 
Moreover let
$$f(z)=\sum_{m=1}^{\infty}c_f(m){\bf e}(mz)$$
be the  primitive form in ${\textfrak S}_{2k-n}(SL_2({\bf Z}))$ corresponding to $g$ via  the Shimura correspondence. 
 For a Dirichlet character $\chi,$ we then define the L-function $L(s,f,\chi)$ of $f$ twisted by $\chi$ by
$$L(s,f,\chi)=\prod_{p}\{(1-\chi(p)\beta_p p^{k-n/2-1/2-s})(1-\chi(p)\beta_p^{-1} p^{k-n/2-1/2-s})\}^{-1},$$
where  $\beta_p$ is  a non-zero complex number such that
$\beta_p+\beta_p^{-1}=p^{-k+n/2+1/2}c_f(p).$ We simply write $L(s,f)$ as $L(s,f,\chi)$ if $\chi$ is the principal character.  
We also define the adjoint $L$-function $L(s,f,{\rm Ad})$ of $f$ by
$$L(s,f,{\rm Ad})=\prod_{p}\{(1-\beta_p^2 p^{-s})(1-\beta_p^{-2} p^{-s})(1-p^{-s})\}^{-1}.$$
We note that $L(s,f,{\rm Ad})$ coincides with $L(s,f,\underline{{\rm St}}).$
Now we put  
 $$c_{I_n(g)}(T)=c_g(|{\textfrak d}_T|)  \prod_p (p^{k-n/2-1/2}\beta_p)^{\nu_p({\textfrak f}_T)}F_p(T,p^{-(n+1)/2}\beta_p^{-1}).$$
We note that $c_{I_n(g)}(T)$ does not depend on the choice of $\beta_p.$ 
Define a Fourier series $I_n(g)(Z)$ by 
$$I_n(g)(Z)= \sum_{T \in {\mathcal H}_n({\bf Z})_{> 0}} c_{I_n(g)}(T){\bf e}({\rm tr}(TZ)).$$ 
In [Ik1] Ikeda showed that  $I_n(g)(Z)$ is a cusp form of weight $k$ with respect to $\varGamma^{(n)}$ and a Hecke eigenform for ${\bf L}_n^{\circ}$ such that 
$$L(s,I_n(g),\underline{\rm St})=\zeta(s)\prod_{i=1}^n L(s+k-i,f).$$
This was first conjecture by Duke and Imamoglu. Thus we call $I_n(g)$ the Duke-Imamoglu-Ikeda lift of $g$  (or of $f$). We note that we have ${\bf Q}(g)={\bf Q}(I_n(g))={\bf Q}(f).$ Furthermore, we have ${\textfrak I}_{g}={\textfrak I}_{I_n(g)},$ where ${\textfrak I}_{g}$ is the ${\textfrak O}_{{\bf Q}(f)}$-module generated by all the Fourier coefficients of $g.$ 

Now to consider a congurence between Duke-Imamoglu-Ikeda  lifts and non-Duke-Imamoglu-Ikeda  lifts, first we prove the following:

\bigskip

{\bf Proposition 4.1} {\it $I_n(g)$ is a Hecke eigenform.}

\bigskip

We note that  Ikeda proved in [Ik1] that $I_n(g)$ is a Hecke eigenform for  ${\bf L}_n^{\circ}$ but has not proved  that it is a Hecke eigenform for ${\bf L}_n.$\footnote[1]{This was  pointed to us by B. Heim (see  [He].) } 
We also note that an explicit form of the spinor L-function of $I_n(g)$ was obtained by Murakawa [Mu] and Schmidt [Sch] assuming that $I_n(g)$ is a Hecke eigenform.

\bigskip

{\bf Proof of Proposition 4.1.} We have only to prove that $I_n(g)$ is an eigenfunction of $T(p)$ for any prime $p.$ The proof may be more or less well konwn, but for the covenience of the readers we here give the proof.  For a modular  form 
$$F(Z)=\sum_B c_F(B){\bf e}({\rm tr}(BZ)),$$
let $c_F^{(p)}(B)$ be the $B$-th Fourier coefficient of $F|T(p).$ Then for any positive definite matrix $B$ we have
$$c_F^{(p)}(B)=p^{nk-n(n+1)/2}\sum_{d_1|d_2|\cdots|d_n|p} d_1^nd_2^{n-1}\cdots d_n$$
$$ \times \sum_{D \in \Lambda_n (d_1 \bot \cdots d_n) \Lambda_n}\det D^{-k} c_F(p^{-1}A[^tD]),$$
where $\Lambda_n=GL_n({\bf Z}).$ 

Now let $E_{n,k}(Z)$ be the Siegel Eisenstein series of degree $n$ and of weight $k$ defined by 
$$E_{n,k}(Z)=\sum_{\gamma \in \varGamma_{\infty}^{(n)} \backslash \varGamma^{(n)}} j(\gamma,Z)^{-k}.$$
For $k \ge n+1,$ the Siegel Eisenstein series $E_{n,k}(Z)$ is a holomorphic modular form of weight $k$ with respect to $\varGamma^{(n)}.$ Furthermore, $E_{n,k}(Z)$ is a Hecke eigenform and in particular we have
$$E_{n,k}|T(p)(Z)=h_{n,p}(p^{k})E_{n,k}(Z),$$
where 
$$h_{n,p}(X)=1+\sum_{r=1}^n \sum_{1 \le i_1 <\cdots <i_r \le n}p^{-\sum_{j=1}^r i_j}X^r.$$
Let $c_{n,k}(B)$ be the $B$-th Fourier coefficient of $E_{n,k}(Z).$ Then we have$$h_{n,p}(p^{k})c_{n,k}(B)=p^{nk-n(n+1)/2}\sum_{d_1|d_2|\cdots|d_n|p} d_1^nd_2^{n-1}\cdots d_n$$
$$ \times \sum_{D \in \Lambda_n (d_1 \bot \cdots d_n) \Lambda_n}\det D^{-k} c_{n,k}(p^{-1}B[^tD]).$$
Let $B$ be positive definite. Then we have
$$c_{n,k}(B)=a_{n,k}(\det 2B)^{k-(n+1)/2}L(k-n/2,\chi_B)\prod_q F_q(B,p^{-k}),$$where $a_{n,k}$ is a non-zero constant depending only on $n$ and $k.$ We note that we have
$$F_q(p^{-1}B[^tD],X)=F_q(B,X)$$
for any $D \in \Lambda_n (d_1 \bot \cdots d_n) \Lambda_n$ with $d_1|\cdots|d_n|p$ if $q \not=p.$ 
Thus we have
$$h_{n,p}(p^{k})F_p(B,p^{-k})=\sum_{e_1 \le e_2 \le \cdots \le e_n \le 1} p^{ne_1+(n-1)e_2+ \cdots + e_n} p^{(e_1+\cdots+e_n)(k-n-1)}$$
$$ \times \sum_{D \in \Lambda_n \backslash \Lambda_n (p^{e_1} \bot \cdots p^{e_n}) \Lambda_n} F_p(p^{-1}B[^tD],p^{-k}).$$
The both-hand sides of the above are polynomilals in $p^{k}$ and the equality holds for infinitely many $k.$ Thus we have
$$h_{n,p}(X^{-1})F_p(B,X)=\sum_{e_1 \le e_2 \le \cdots \le e_n \le 1} p^{ne_1+(n-1)e_2+ \cdots + e_n} (X^{-1}p^{-n-1})^{(e_1+\cdots+e_n)}$$
$$ \times \sum_{D \in \Lambda_n \backslash \Lambda_n (p^{e_1} \bot \cdots p^{e_n}) \Lambda_n} F_p(p^{-1}B[^tD],X)$$
as polynomials in $X$ and $X^{-1}.$ Thus we have
$$ (p^{k-(n+1)/2}X)^{n/2} h_{n,p}(p^{(n+1)/2}X^{-1})  (p^{k-(n+1)/2}X^{-1})^{\nu_p({\textfrak f}_B)} F_p(B,p^{-(n+1)/2}X)$$
$$=p^{nk-n(n+1)/2}\sum_{e_1 \le e_2 \le \cdots \le e_n \le 1} p^{ne_1+(n-1)e_2+ \cdots + e_n}$$
$$ \times \sum_{D \in \Lambda_n (p^{e_1} \bot \cdots p^{e_n}) \Lambda_n}\det D^{-k} (p^{k-(n+1)/2}X^{-1})^{\nu_p({\textfrak f}_{p^{-1}B[^tD] )} }F_p(p^{-1}B[^tD],p^{-(n+1)/2}X).$$
We recall that we have
$$c_{I_n(g)}(B)=c_{g}(|{\textfrak d}_B|) {\textfrak f}_B^{k-(n+1)/2} \prod_q (\beta_q)^{\nu_q({\textfrak f}_B)}F_q(B,q^{-(n+1)/2}\beta_q^{-1}),$$
where for each prime number $q$ $\beta_q$ is the non-complex number defined before. 
We also note that $c_{g}(|{\textfrak d}_{p^{-1}B[^tD]}|)=c_{g}(|{\textfrak d}_B|)$ for any $D.$ Thus we have
$$ (p^{k-(n+1)/2}\alpha_p^{-1})^{n/2} h_{n,p}(p^{(n+1)/2}\alpha_p)  c_{I_n(g)}(B)$$
$$=p^{nk-n(n+1)/2}\sum_{d_1|d_2|\cdots|d_n|p} d_1^nd_2^{n-1}\cdots d_n\sum_{D \in \Lambda_n (d_1 \bot \cdots d_n) \Lambda_n}\det D^{-k} c_{I_n(g)}(p^{-1}B[^tD]).$$
This proves the assertion.

\bigskip

Let $\{f_1,....,f_d \}$ be a basis of  ${\textfrak S}_{2k-n}(\varGamma^{(1)})$ consisting of primitive forms. Let $K$ be an algebraic number field containing ${\bf Q}(f_1)\cdots {\bf Q}(f_d),$ and $A={\textfrak O}_K.$  
To formulate our conjecture exactly, we introduce the Eichler-Shimura periods as follows (cf. Hida [Hi3].)   
Let $f$ be a primitive form in ${\textfrak S}_{2k-n}(\varGamma^{(1)}).$ 
Let $\textfrak P$ be a prime ideal in $K.$ Let $A_{\textfrak P}$ be a valuation ring in $K$ corresponding to $\textfrak P.$ Assume that the residual characteristic of $A_{\textfrak P}$ is greater than or equal to $5.$ Let $L(2k-n-2,A_{\textfrak P})$ be the module of homogeneous polynomials of degree $2k-n-2$ in the variables $X,Y$ with coefficients in $A_{\textfrak P}.$ We define the action of $M_2({\bf Z}) \cap GL_2({\bf Q})$  on $L(2k-n-2,A_{\textfrak P})$ via 
$$\gamma \cdot P(X,Y)=P({}^t(X,Y)(\gamma)^{\iota}),$$
where $\gamma^{\iota}=(\det \gamma)\gamma^{-1}.$ 
Let $H_P^1(\varGamma^{(1)} ,L(2k-n-2,A_{\textfrak P}))$ be the parabolic cohomology group of $\varGamma^{(1)}$ with values in $L(2k-n-2,A_{\textfrak P}).$ Fix a point $z_0 \in {\bf H}_1.$ Let $g \in {\textfrak S}_{2k-n}(\varGamma^{(1)})$ or $g \in \overline{{\textfrak S}_{2k-n}(\varGamma^{(1)})}.$ We then define the differential $\omega(g)$ as 
$$\omega(g)(z)=\left\{\begin{array} {ll}
2\pi i g(z) (X-zY)^n dz & \ {\rm if} \ g \in {\textfrak S}_{2k-n}(\varGamma^{(1)}) \\
2\pi \sqrt{-1} g(z) (X-\bar zY)^n dz & \ {\rm if} \ g \in \overline{{\textfrak S}_{2k-n}(\varGamma^{(1)})},
\end{array} \right.$$
 and define the cohomology class $\delta(g)$ of the  1-cocycle of $\varGamma^{(1)}.$ as
$$\gamma \in \varGamma^{(1)}  \longrightarrow \int_{z_0}^{\gamma(z_0)} \omega(g).$$ The mapping $g \longrightarrow \delta(g)$ induces the isomorphism  
$$\delta:{\textfrak S}_{2k-n}(\varGamma^{(1)}) \oplus \overline{{\textfrak S}_{2k-n}(\varGamma^{(1)})} \longrightarrow H_P^1(\varGamma^{(1)} ,L(2k-n-2,{\bf C})),$$
which is called the Eichler-Shimura isomorphism. 
We can define the action of Hecke algebra ${\bf L}_1'$ on $H_P^1(\varGamma^{(1)}, L(2k-n-2,A_{\textfrak P}))$ in a natural manner. Furthermore, we can  define the action $F_{\infty}$ on $H_P^1(\varGamma^{(1)}, L(2k-n-2,A_{\textfrak P}))$ as
$$F_{\infty}(\delta(g)(z))=\mattwo(-1;0;0;1)\delta(g)(-\bar z),$$
and this action commutes with the Hecke action.   For a primitive form $f$ and $j=\pm 1,$ we define the subspace $H_P^1(\varGamma^{(1)},L(2k-n-2,A_{\textfrak P}))[f,j]$ of $H_P^1(\varGamma^{(1)}, L(2k-n-2,A_{\textfrak P}))$ as
$$H_P^1(\varGamma^{(1)},L(2k-n-2,A_{\textfrak P}))[f,j]$$
$$=\{ x \in H_P^1(\varGamma^{(1)} , L(2k-n-2,A_{\textfrak P})) \ ; \ x|T=\lambda_f(T)x \ {\rm for } \ T \in {\bf L}_1, \ {\rm and} \ F_{\infty}(x)=jx \}.$$
Since $A_{\textfrak P}$ is a principal ideal domain, $H_P^1(\varGamma^{(1)} ,L(2k-n-2,A_{\textfrak P}))[f,j] $ is a free module of rank one over $A_{\textfrak P}.$ For each $j=\pm 1$ take a basis $\eta(f,j,A_{\textfrak P})$ of $H_P^1(\varGamma^{(1)},(2k-n-2,A_{\textfrak P}))[f,j] $ and define a complex number $\Omega(f,j;A_{\textfrak P})$ by
$$(\delta(f)+jF_{\infty}(\delta(f)))/2=\Omega(f,j;A_{\textfrak P})\eta(f,j;A_{\textfrak P}).$$
This $\Omega(f,j;A_{\textfrak P})$ is uniquely determined up to constant multiple of units in $A_{\textfrak P}.$ We call $\Omega(f,+;A_{\textfrak P})$ and $\Omega(f,-;A_{\textfrak P})$ the Eichler-Shimura periods. 
For $j=\pm, 1 \le l \le 2k-n-1,$  and a Dirichlet character $\chi$ such that $\chi(-1)=j (-1)^{l-1},$ put 
$${\bf L}(l,f,\chi)={\bf L}(l,f,\chi;A_{\textfrak P})={\Gamma(l) L(l,f,\chi) \over \tau(\chi)(2\pi \sqrt{-1})^{l}\Omega(f,j;A_{\textfrak P})},$$
where $\tau(\chi)$ is the Gauss  sum of $\chi.$ 
In particular, put ${\bf L}(l,f;A_{\textfrak P})={\bf L}(l,f,\chi;\textfrak P)$ if $\chi$ is the principal character. Furthermore, put 
$$\Gamma_{\bf C}(s)=2(2\pi)^{-s}\Gamma(s),$$ and
$${\bf L}(s,f,{\rm Ad})={\Gamma_{\bf C}(s)\Gamma_{\bf C}(s+2k-n-1)L(s,f, {\rm Ad}) \over \langle f, f \rangle}.$$
It is well-known that ${\bf L}(l,f,\chi)$ belongs to the field $K(\chi)$ generated over $K$ by all the values of $\chi,$  and ${\bf L}(l,f,\underline {\rm St})$ belongs to ${\bf Q}(f)$ (cf. [Bo].) Let $I_n(g)$ be the Duke-Imamoglu-Ikeda  lift of $f.$ Let ${\textfrak S}_k(\varGamma^{(n)})^*$ be the subspace of ${\textfrak S}_k(\varGamma^{(n)})$ generated by all the Duke-Imamoglu-Ikeda  lifts $I(g)^n$ of  primitive forms $g \in {\textfrak S}_{2k-n}(\varGamma^{(1)}).$ We remark that ${\textfrak S}_k(\varGamma^{(2)})^*$ is the Maass subspace of ${\textfrak S}_k(\varGamma^{(2)}).$ 

\bigskip

{\bf Conjecture A.} {\it Let  $g$ be a Heke eigenform belonging to the Kohnen plus subspace ${\textfrak S}_{k-n/2+1/2}^+(\varGamma_0(4)),$ and  $f$ the primitive form in ${\textfrak S}_{2k-n}(\varGamma^{(1)})$ corresponding to $g$ via the Shimura correspondence. Moreover let $K$ be the field as above. Assume that $k >n.$ Let  ${\textfrak P}$ be a prime ideal of $K$ not dividing $(2k-1)!.$  Then  ${\textfrak P}$ is a congruence prime of$I_n(g)$ with respect to $({\textfrak S}_k(\varGamma^{(n)})^*)^{\bot} $ if ${\textfrak P}$ divides  \\
${\bf L}(k,f)\prod_{i=1}^{n/2-1} {\bf L}(2i+1,f,{\rm Ad}) .$}

\bigskip

{\bf Remark.} This is an analogue of the Doi-Hida-Ishii conjecture concerning the  congruence  primes of the Doi-Naganuma lifting [D-H-I]. (See also [Ka1].) We also note that this type of conjecture has been proposed by Harder [Ha] in the case of vector valued Siegel modular forms. 

\bigskip 

Now to explain why our conjecture is reasonable, we refer to Ikeda's conjecture on the Petersson inner product of the Duke-Imamoglu-Ikeda  lift. Let $g$ and $f$  be as above. Put
$$\tilde \xi(s)=\Gamma_{\bf C}(s)\zeta(s).$$

\bigskip

{\bf Theorem 4.2.} (Katsurada and Kawamura [K-K]) {\it Under the above notation and the assumption we have 
$$\tilde \xi(n)\Gamma_{\bf C}(k)L(k,f)\prod_{i=1}^{n/2-1} {\bf L}(2i-1,f,{\rm Ad})\tilde \xi(2i) 
=2^{\alpha} { \langle I_n(g) f,I_n(g) \rangle  \over \langle f,f \rangle^{n/2-1} \langle g,g \rangle },$$
where $\alpha$ is an integer depending only on $n$ and $k.$ }

\bigskip

We note that the above theorem was conjectured by Ikeda [Ik2] under more general setting.  We note that the theorem has been proved by Kohnen and Skoruppa [K-S] in case $n=2.$ 

\bigskip

{\bf Proposition 4.3} {\it Let the notation and the assumption be as above. For a fundamental discriminant $D$ such that $(-1)^{n/2}D >0$ let $\chi_D$ be the Kronecker character correponding to $D.$ Then we have 
$${c_g(|D|)^2 \langle  f,\, f \rangle ^{n/2} \over \langle I_n(g),\, I_n(g) \rangle}= {2^{a_{n,k}}(-1)^{b_{n,k}}  |D|^{k-n/2}{\bf L}(k-n/2,f,\chi_D) \over {\bf L}(k,\, f)  \widetilde{\xi}(n) \displaystyle\prod_{i=1}^{n/2-1} {\bf L}(2i+1,\, f,\, {\rm Ad})  \widetilde{\xi}(2i)}$$
with some integers $a_{n,k}$ and $b_{n,k}$ depending only on $n$ and $k.$ 
}

\bigskip

{\it Proof.} By the result in Kohnen-Zagier[K-Z], for any fundamental discriminant $D$ such that $(-1)^{n/2}D >0$ we have
$${c_g(|D|)^2 \over \langle g, g \rangle} ={2^{k-n/2-1}|D|^{k-n/2-1/2} \Gamma_{\bf C}(k-n/2)L(k-n/2,f,\chi_D) \over \langle f, f \rangle}.$$
We note that $\tau(\chi_D)$ is $\sqrt {D}$ or $\sqrt{-1}\sqrt{D}$ according as $n \equiv 0 \ {\rm mod} \ 4,$ or $n \equiv 2 \ {\rm mod} \ 4.$ This completes the proof.

\bigskip

{\bf Lemma 4.4.} {\it Let $f$ be as above.

{\rm (1)} Let  ${\bf r}_1$ be an element of ${\bf L}_n'$ in Proposition 2.3. Then we have
$$\lambda_{I_n(g)}({\bf r}_1)=p^{(n-1)k-n(n+1)/2}c_f(p)\sum_{i=1}^n p^{i}.$$

{\rm (2)} Let $n=2.$ Then we have
$$\lambda_{I_2(g)}(T(p))=c_f(p)+p^{2k-n-1}+p^{2k-n-2}.$$

}

\bigskip

{\bf Lemma 4.5.} {\it  Let $d$ be a fundamental discriminant such that $(-1)^{n/2}d >0.$ 

\noindent
{\rm (1)} Assume that $d \not=1.$  Then there exists a positive definite half integral matrix $A$ of degree $n$ such that $(-1)^{n/2}\det (2A)=d.$

\noindent
{\rm (2)} Assume  $n \equiv 0 \ {\rm mod} \ 8.$ Then there exsits a  positive definite half integral matrix $A$ of degree $n$ such that $\det (2A)=1.$

\noindent
{\rm (3)} Assume that $n \equiv 4 \ {\rm mod} \ 8.$ Then for any prime number $q$ there exsits a  positive definite half integral matrix $A$ of degree $n$ such that $det (2A)=q^2.$ 

 }

\bigskip

{\it Proof.} (1) For two elements $a,b$ of ${\bf Q}_p^{\times},$ let $(a,b)_p$ denote the Hilbert symbol. For a non-degenerate symmetric matrix $A$ with entries in ${\bf Q}_p$ let $h_p(A)$ be the Hasse invariant of $A.$ (For the definition of the Hasse invariant, see, for example, [I-S].) First let $n \equiv 2 \ {\rm mod} \ 4$ and $d=-4.$ Take a family $\{A_p \}_p$ of half integral matrices over ${\bf Z}_p$ of degree $n$ such that
$A_p=  1_{n}$ if $p \not= 2,$ and $A_2=(-1)^{(n-2)/4}1_2 \bot H_{n/2-1},$ 
where 
$H_r = \overbrace{H \bot ...\bot H}^r$ with 
$H =\left(\begin{array}{cc}
             0        & 1/2  \\
             1/2        & 0 
                                   \end{array}\right).$ Then we have $\det A = 2^{2-n}  \in {{\bf Q}_p}^{\times} /({{\bf Q}_p}^{\times})^{2}$ for any $p,$ and  $h_p(A)=1$ for any $p.$ 
Thus by [I-S, Proposition 2.1], there exists an element $A$ of ${{\mathcal L}_{n,2}}_{>0}$ such that $A \sim A_p$ for any $p.$ In paricular we have $(-1)^{n/2} \det (2A)=-4.$ Next let $d=(-1)^{n/2}8.$ We take $A_p=(-1)^{n/2}2 \bot 1_{n-1}$ if $p \not=2.$ We can take $\xi \in {\bf Z}_2^*$ such that $(2,\xi)_2=(-1)^{(n-2)(n+4)/8},$ and put $A_2=2\xi \bot (-\xi) \bot H_{n/2-1}.$  Then we have
$\det A = (-1)^{n/2}2^{3-n}  \in {{\bf Q}_p}^{\times} /({{\bf Q}_p}^{\times})^{2}$ for any $p,$ and  $h_p(A)=1$ for any $p.$ Thus again by [I-S, Proposition 2.1], we prove the assertion for this case. Finally assume that $d$ contains a odd prime factor $q.$ For $p \not =q$ we take a matrix $A_p$ so that $\det A_p=2^{-n}d \in {{\bf Q}_p}^{\times} /({{\bf Q}_p}^{\times})^{2}.$ Then for almost all $p$ we have $h_p(A_p)=1.$ We take $\xi \in {\bf Z}_q^*$ such that $(q,-\xi)_q=
\prod_{p \not=q} h_p(A_p),$ and put $A_q=\xi d \bot \xi \bot 1_{n-2}.$ Then we have $2^{-n}d \det A_q  \in {{\bf Q}_q}^{\times} /({{\bf Q}_q}^{\times})^{2},$ and $h_q(A_q) \prod_{p \not=q}h_p(A_p)=1.$ Thus again by [I-S, Proposition 2.1], we prove the assertion for this case.

(2) It is well known that there exists a positive definite half-integral matrix $E_8$ of degree $8$ such that $\det (2E_8)=1.$ Thus $A=\overbrace{E_8 \bot \cdots \bot E_8}^{n/8}$ satisfies the required condition.

(3) Let $q \not=2.$ Then, take a family $\{A_p \}$ of half-inegral matrices over ${\bf Z}_p$ of degree $n$ such that  $A_q \sim_{{\bf Z}_q} q \bot (-q\xi) \bot (-\xi) \bot 1_{n-3}$ with $({\xi \over q})=-1, A_2=H_{n/2},$ and $A_p=1_n$ for $p \not= q,2.$ Then by the same argumemnt as in (1) we can show that there exits a  positive definite half integral matrix $A$ of degree $n$ such that $\det (2A)=q^2$  such that $A \sim_{{\bf Z}_p} A_p$ for any $p.$ Let $q=2.$ Then the matrix 
$A'=\Biggl(\begin{smallmatrix}
1    & 0   & 0  & 1/2\\
0    & 1   & 0  & 1/2\\
0    & 0   & 1  & 1/2 \\
1/2  & 1/2 & 1/2& 1
\end{smallmatrix}\Biggr ) $
 is a positive definite and $\det (2A')=4.$ Thus the matrix $A' \bot \overbrace{E_8 \bot \cdots \bot E_8}^{(n-4)/8}$ satisfies the required condition.

\bigskip

\bigskip

{\bf Proposition 4.6.}{\it Let $k$ and $n$ be positive even integer. Let $d$ be a fundamental discriminant. Let $f$ be a primitive form in ${\textfrak S}_{2k-n}(\varGamma^{(1)}).$  Let $\textfrak P$ be a prime ideal in $K.$ Then there exists a positive definite half integral matrices $A$ of degree $n$ such that $c_{I_n(g)}(A)=c_{g}(|d|)l$ with an integer $l$ not divisible by $\textfrak P.$
}

\bigskip

{\it Proof.} First assume that $d \not=1,$ or $n \not \equiv 4 \ {\rm mod } \ 8.$ (1) By (1) and (2) of Lemma 4.5, there exists a matrix $A$ such that 
${\textfrak d}_A=d.$ Thus we have $c_{I_n(g)}(A)=c_{g}(|d|).$ This proves the assertion. 

Next assume that  $n \equiv 4 \ {\rm mod} \ 8$ and that $d=1.$ We show that there exists a prime number $q$ satisfying the following condition: 

(**) $c_f(q)+q^{k-n/2-1}(-q-1)$ is not divisible by $\textfrak P.$\footnote[2]{The proof of this fact  was  suggested by S. Yasuda and T. Yamauchi.}   

Assume that $c_f(q)+q^{k-n/2-1}(-q-1)$ is  divisible by $\textfrak P$ for any prime number $q.$ Let $p$ be a prime number divisible by ${\textfrak P}.$ Fix an imbedding $\iota_p:\bar {\bf Q} \longrightarrow \overline {{\bf Q}_p},$ and let $\rho_{f,p}:Gal(\bar {\bf Q}/{\bf Q}) \longrightarrow GL_2(\overline{{\bf Q}}_p)$ be the Galois representation attached to $f.$ Then by Chebotarev density theorem, the semi-simplification $\overline{\rho}_{f,p}^{ss}$ of $\overline{\rho}_{f,p}$ can be expressed as
 $$\overline{\rho}_{f,p}^{ss}=\overline{\chi_p}^{k-n/2} \oplus \overline{\chi_p}^{k-n/2-1}$$ 
 with $\overline{\chi_p}$ the $p$-adic mod $p$ cyclotomic character. On the other hand, by the Fontaine-Messing [Fo-Me] and Fontaine-Laffaille [Fo-La], $\overline{\rho}_{f,p}^{ss}|I_p$ should be $\overline{\chi_p}^{2k-n-1} \oplus 1$ or ${\omega}_2^{2k-n-1} \oplus {\omega}_2^{p(2k-n-1)}$ with ${\omega}_2$ the fundamental character of level 2, where $I_p$ denotes the inertia group of $p$ in $Gal(\bar {\bf Q}/{\bf Q}).$ This is impossible because $k >2.$ Now for a prime number $q$ satisfying the condition (**), take a positive definite matrix $A$ in (3) of Lemma 4.5. 
 Then
 $$c_{I_n(g)}(A)=c_g(1) q^{k-(n+1)/2} \beta_qF_q(A,q^{-(n+1)/2}\beta_q^{-1}).$$ By [Ka1], we have
$$F_q(B,X)=1- Xq^{(n-2)/2}(q^2+ q)+q^3(Xq^{(n-2)/2})^2.$$
Thus we have
$$c_{I_n(g)}(A)=c_g(1)(c_f(q)+q^{k-n/2-1}(-q-1)).$$
Thus the assertion holds.

\bigskip

\bigskip

{\bf Theorem 4.7.} {\it Let $k \ge 2n+4.$ Let $K$ and $f$ be as above, and  ${\textfrak P}$ a prime ideal of $K.$ Furthermore assume that \\
{\rm (1)} ${\textfrak P}$ divides ${\bf L}(k,f) \prod_{i=1}^{n/2-1} {\bf L}(2i+1,f, {\rm Ad}). $ \\
{\rm (2)} ${\textfrak P}$ does not divide 
$$\tilde \xi(2m)\prod_{i=1}^{n} {\bf L}(2m+k-i,f)
{\bf L}(k-n/2,f,\chi_D)D(2k-1)!$$ 
for some integer $n/2 +1  \le m \le k/2-n/2-1,$  and for some fundamental discriminant $D$ such that $(-1)^{n/2}D>0.$  \\
Then ${\textfrak P}$ is a congruence prime of $I_n(g)$ with respect to ${\bf C} {I_n(g)}^{\bot}.$ Furthermore assume that the following  condition hold:\\

{\rm (3)} ${\textfrak P}$ does not divide 
 $$C_{k,n}{\langle f, f \rangle \over \Omega(f,+,A_{\textfrak P})\Omega(f,-,A_{\textfrak P})},$$
 where $C_{k,n}=1$ or $\prod_{q \le (2k-n)/12} (1+q+\cdots+q^{n-1})$ according as $n=2$ or not.\\
 
Then ${\textfrak P}$ is a congruence prime of $I_n(g)$ with respect to $({\textfrak S}_k(\varGamma^{(n)})^*)^{\bot}.$}

\bigskip

{\it Proof.}   Let $\textfrak P$ be a prime ideal satisfying the condition (1) and (2). For the $D$ above, take a matrix $A \in {\mathcal H}_n({\bf Z})_{>0}$ so that $c_{I_n(g)}(A)=c_{g}(|D|)l$ with an integer $l$ not divisible by ${\textfrak P}.$ Then by Proposition 4.3, we have 
$$\Lambda(2m,I_n(g),\underline {\rm St})|c_{I_n(g)}(A)|^2=\Lambda(2m,I_n(g),\underline {\rm St})|c_{g}(|D|)|^2l^2$$
$$=\epsilon_{k,m}{\prod_{i=1}^n {\bf L}(2m+k-i,f) |D|^{k-n/2}{\bf L}(k-n/2, f,\chi_D) \over {\bf L}(k,f) \tilde \xi(n)\prod_{i=1}^{n/2-1}{\bf L}(2i+1,f, {\rm Ad})\tilde \xi(2i)} $$
$$ \times ({ \displaystyle \Omega(f,+;\textfrak P)\Omega(f,-;A_{\textfrak P}) \over \displaystyle  \langle f,f \rangle  })^{n/2},$$
where $\epsilon_{k,m}$ is a rational number whose numerator is not divided by ${\textfrak P}.$ 
We note that ${ \displaystyle \langle f,f \rangle  \over \displaystyle \Omega(f,+;A_{\textfrak P})\Omega(f,-;A_{\textfrak P})}$ is ${\textfrak P}$-integral.
Thus by assumptions (1) and (2), 
$\textfrak P$ divides $(\Lambda(2m,I_n(g),\underline {\rm St})c_{I_n(g)}(A)^2)^{-1},$ and thus it divides $(\Lambda(2m,I_n(g),\underline {\rm St}){\textfrak I}_{I_n(g)}^2)^{-1}.$ We note that $I_n(g)$ satisfies the assumption in Theorem 3.1. Thus by Theorem 3.1, there exits a Hecke eigenform $G \in {\bf C}(I_n(g))^{\bot}$ such that 
$$\lambda_G(T) \equiv \lambda_{I_n(g)}(T) \ {\rm mod} \ \textfrak P$$
 for any $T \in {\bf L}'_n.$  Assume that we have $G =I_n(g')$ with some Hecke eigenform $g' \not=g$ in ${\textfrak S}_{k-n/2+1/2}^+(\varGamma_0(4)),$ and let $f'(z)=\sum_{m=1}^{\infty} c_{f'}(m) {\bf e}(mz)$  be the primitive form in ${\textfrak S}_{2k-n}(\varGamma^{(1)})$ corresponding to $g'$ via the Shimura correpondence. Let $n=2.$ Then by (2) of Lemma 4.4, $\textfrak P$ is also a congruence prime of $f.$ Let $n \ge 4.$ Then by (1) of Lemma 4.4, we have
 $$(p^{n-1}+\cdots +p+1)c_f(p) \equiv (p^{n-1}+\cdots +p+1)c_{f'}(p) \ {\rm mod} \ \textfrak P$$
 for any prime number $p$ not divisible by $\textfrak P.$ By assumption (3), in particular, for any $p \le (2k-n)/12,$ we have
 $$c_f(p) \equiv c_{f'}(p) \ {\rm mod} \ \textfrak P.$$
 Thus by Sturm [Stur], $\textfrak P$ is also a congruence prime of $f.$ Thus by [Hi2] and [Ri2], $\textfrak P$ divides ${ \displaystyle \langle f,f \rangle  \over \displaystyle \Omega(f,+;A_{\textfrak P})\Omega(f,-;A_{\textfrak P})},$ which contradicts the assumption (3). Thus ${\textfrak P}$ is a congruence prime of $I_n(g)$ with respect to $({\textfrak S}_{k}(\varGamma^{(n)})^*)^{\bot}.$

\hskip 8cm $\Box$

\bigskip

{\bf Example} Let $n=4$ and $k=18.$ Then we have 
$\dim \ {\textfrak S}_{18}(\varGamma^{(4)}) \approx 16$  
(cf. Poor and Yuen[P-Y]) and $\dim \ {\textfrak S}_{18}(\varGamma^{(4)})^*=\dim \ {\textfrak S}_{17}^+(\varGamma_0(4))= \dim \ {\textfrak S}_{32}(\varGamma^{(1)})=2.$ 
Take a basis $g_1,g_2$ consisting of Hecke eigenforms in ${\textfrak S}_{17}^+(\varGamma_0(4)),$  and
$f_i$ be the primitive form in ${\textfrak S}_{32}(\varGamma^{(1)})$ corresponding to $g_i$ under the Shimura correspondence. Put $g=g_1$ and $f=f_1.$ Then we have $[{\bf Q}(f):{\bf Q}]=2.$ Now consider the prime number $211.$ Then it is decomposed into two prime ideals as $211={\textfrak P}{\textfrak P}'$ in ${\bf Q}(f).$ Then we have 
$$N_{{\bf Q}(f)/{\bf Q}}({\bf L}(18,f))=2^7 \cdot 3^2 \cdot 5^2 \cdot 7^2 \cdot 11 \cdot 13 \cdot 211,$$
$$N_{{\bf Q}(f)/{\bf Q}}(\prod_{i=1}^4 {\bf L}(24-i,f))=2^{19} \cdot 3^{13} \cdot 5^5 \cdot 7^8 \cdot 11^2 \cdot 13^5 \cdot 17^5 \cdot 19^3 \cdot 23 \cdot 503 \cdot 1307 \cdot 14243,$$  
$$\tilde \xi(6)=2^{-2} \cdot 3^{-2} \cdot 7^{-1}$$
and
$$N_{{\bf Q}(f)/{\bf Q}}({\bf L}(16,f,\chi_1)=2^5 \cdot 3^2 \cdot 5^3 \cdot 7^2 \cdot 11 \cdot 13^2.$$
(cf. Stein [Ste].)  Thus by Theorem 4.7, ${\textfrak P}$ or ${\textfrak P}'$ is a congruence prime of $I_n(g)$ with respect to $({\bf C}I_n(g))^{\bot}.$ Furthermore, by a direct computation we see neither ${\textfrak P}$ nor ${\textfrak P}'$ is a congruence prime of $I_n(g)$ with respect to ${\bf C}I_n(g_2).$ This implies that ${\textfrak P}$ or ${\textfrak P}'$ is a congruence prime of $I_n(g)$ with respect to  ${{\textfrak S}_{18}(\varGamma^{(4)})^*}^{\bot}.$

\bigskip

{\bf References}

\noindent
[A] A. N. Andrianov, Quadratic forms and Hecke operators, Springer, 1987.

\noindent
[Bo] S. B{\"o}cherer, {\"U}ber die Fourier-Jacobi-Entwicklung Siegelscher 
Eisensteinreihen II, Math. Z. 189(1985), 81-110.

\noindent
[Br] J. Brown, Saito-Kurokawa lifts and applications to the Bloch-Kato conjecture,  Compos. Math. 143 (2007), 290--322.

\noindent
[D-H-I] K. Doi, H. Hida, and H. Ishii, Discriminant of Hecke fields and twisted adjoint L-values for GL(2), Invent. Math. 134(1998), 547-577.

\noindent
[Fo-La] J. M. Fontaine and G. Laffaille, Construction de repr\'esentation $p$ adiques, Ann. Sci. Math. \'Ecole Norm. Sup. 15(1982), 547-608.

\noindent
[Fo-Me] J. M. Fontaine and W. Messing, $p$-adic periods and $p$-adic \'etale cohomology, Contemp. Math. 67, 179-207.

\noindent
[Ha] G. Harder, A congruence between a Siegel modular form and an elliptic modular form, Preprint 2003.

\noindent
[He] B. Heim, Miyawaki's $F_{12}$ spinor L-function conjecture, \\
arXiv:0712.1286v1[math.NT] 8.12.2007.

\noindent
[Hi1]  H. Hida, Congruences of cusp forms and special values of their zeta functions,  Invent Math. 63(1981), 225-261.

\noindent
[Hi2]  H. Hida, On congruence divisors of cusp forms as factors of the special avlues of their zeta functions, Invent. Math. 64(1981), 221-262

\noindent
[Hi3] H. Hida, Modular forms and Galois cohomology, Cambridge Univ. Press, 2000

\noindent
[Ik1] T. Ikeda, On the lifting of elliptic modular forms to Siegel cusp forms of degree $2n,$ Ann. of Math.@154(2001), 641-681.

\noindent
[Ik2] T, Ikeda, Pullback of lifting of elliptic cusp forms and Miyawaki's conjecture, Duke Math. J. 131 (2006), 469-497.

\noindent
[Ka1] H. Katsurada, Special values of the standard zeta functions for elliptic modular forms, Experiment. Math. 14(2005), 27-45. 

\noindent
[Ka2] H. Katsurada, Congruence of Siegel modular forms  and  special values of their  standard zeta functions, Math. Z. 259 (2008), 97-111. 

\noindent
[Ka3] H. Katsurada, A remark on the normalization of the standard zeta values for Siegel modular forms, Abh. Math. Sem. Univ. Hamburg 80(2010), 37-45.

\noindent
[K-K] \underline{\hspace{22.5mm}}, Ikeda's conjecture on the period of the Duke-Imamoglu-Ikeda  lift, preprint.

\noindent
[K-S] W. Kohnen and N-P. Skoruppa, A certain Dirichlet series attached to Siegel modular forms of degree 2, Invent. Math. 95(1989), 541-558.

\noindent
[K-Z] W. Kohnen and D. Zagier, Values of L-series of modular forms at the center of the critical strip, Invent. Math. 64(1981), 175-198

\noindent
[Mi1] S. Mizumoto, Poles and residues of standard L-functions attached to Siegel modular forms, Math. Ann. 289(1991) 589-612.

\noindent
[Mi2] \underline{\hspace{22.5mm}}, On integrality of Eisenstein liftings, Manuscripta Math. 89(1996), 203-235. \\
Corrections Ibid.307(1997), 169-171.

\noindent
[Mu] K. Murakawa, Relations between symmetric power $L$-functions and spinor $L$-functions attched to Ikeda lifts, Kodai Math. J. 25(2002), 61-71

\noindent
[P-Y] C. Poor and D. Yuen, Private communication (2005).

\noindent
[Ri1] K. Ribet, A modular construction of unramified $p$-extesnions of ${\bf Q}(\nu_p),$ Invent. Math. 34(1976), 151-162.

\noindent
[Ri2] \underline{\hspace{22.5mm}}, Mod $p$ Hecke operators and congruence between modular forms,
 Invent. Math. 71(1983), 193-205.

\noindent
[Sch] R. Schmidt,  On the spin $L$-function of Ikeda's lifts. Comment. Math. Univ. St. Pauli 52 (2003), 1-46.  

\noindent
[Sh] G. Shimura, The special values of the zeta functions associated with cusp forms, Comm. pure appl. Math. 29(1976), 783-804. 

\noindent
[Ste] W. A. Stein, The modular forms data base, \\
http://modular.fas.harvard.edu/index.html

\noindent
[Stur] J. Sturm, Congruence of modular forms, Springer Lect. Notes in Math. 1240(1984) 275-280.

\bigskip

Hidenori Katsurada \\ 
Muroran Institute of Technology \\
27-1 Mizumoto Muroran, 050-8585, Japan \\
e-mail: hidenori@mmm.muroran-it.ac.jp

\end{document}